\documentclass[10pt]{article}
\usepackage{amsfonts, latexsym,amscd,graphicx,amssymb}

\makeatletter \@addtoreset{equation}{section} \makeatother

\def\text#1{\mbox{\rm #1}}
\newcommand{\be}{\begin{equation}}
\newcommand{\ee}{\end{equation}}
\newcommand{\ba}{\begin{eqnarray}}
\newcommand{\ea}{\end{eqnarray}}

\newcommand{\pa}{\partial}
\newcommand{\f}{\frac}

\begin{document}

\title{  Liapunov's direct method for Birkhoffian systems: Applications
to electrical networks}
\author{Delia Ionescu-Kruse,\\
{\small {\it Institute of Mathematics of the Romanian Academy}}\\
{\small {\it P.O. Box 1-764, RO-014700, Bucharest,
 Romania,Delia.Ionescu@imar.ro}}}
 \date{}
 \maketitle
\begin{abstract}
In this paper, the concepts and the direct theorems of stability
in the sense of Liapunov, within the framework of Birkhoffian
dynamical systems on manifolds, are  considered. The Liapunov-type
functions are constructed for linear and nonlinear LC and RLC
electrical networks,  to prove stability under certain conditions.

\textit{Keywords}:  Liapunov stability, geometric methods in
differential equations, Birkhoffian differential systems,
electrical networks.

MSC: 34A26, 34D20, 58A20, 94C
\end{abstract}

\section{Introduction}

During the last years, a far reaching generalization of the
Hamiltonian framework has been developed in a series of papers.
This generalization, which is based on the geometric notion of
generalized Dirac structure (see  Courant \cite{courant} and
Dorfman \cite{dorfman}), gives rise to implicit Hamiltonian
systems  (see, for example, the papers by   Maschke and van der
Schaft \cite{maschke2}, \cite{schaft}). Applications to
 nonholonomic systems  and electrical circuits
 (see Bloch and Crouch \cite{bloch},  Maschke and  van der Schaft
 \cite{maschke2})
 illustrate this theory.
Recently, the notion of implicit Lagrangian system has been
developed by Yoshimura and Marsden \cite{yoshimura}. Nonholonomic
mechanical systems and degenerate Lagrangian systems such as LC
circuits can be systematically formulated in the implicit
Lagrangian context in which Dirac structures are also used.

An alternative approach to the study of dynamical systems
  is the Birkhoffian formalism.
 This is  a global formalism of implicit systems of
second order ordinary differential equations on a manifold. It
applies to a wide class of systems, among them,  nonholonomic
systems,  degenerate systems as well as dissipative systems.
Kobayashi and Oliva developed in \cite{oliva} the framework of
Birkhoffian dynamical systems on manifolds, following Birkhoff's
ideas presented locally in \cite{bir}. The space of configurations
is a smooth $m$-dimensional differentiable connected manifold  and
the covariant character of the Birkhoff generalized forces is
obtained by defining the notion of elementary work, called
Birkhoffian, a special Pfaffian form defined  on the 2-jets
manifold.   The dynamical system associated to this Pfaffian form
is a subset of the 2-jets manifold which defines an implicit
second order ordinary differential system. The notion of
Birkhoffian allows the introduction of the intrinsic concepts of
reciprocity, regularity, affine structure in the accelerations,
conservativeness \cite{oliva}, dissipativeness \cite{eu2}.

The  Birkhoffian formalism in the context of electrical circuits
was discussed  by Ionescu and Scheurle \cite{io} for the case of
LC circuits, and Ionescu \cite{eu2} for the case of RLC circuits.
An LC/RLC circuit, with no assumptions placed on its topology,
will be described by a family of Birkhoffian systems,
parameterized by a finite number of real constants which
correspond to initial values of certain state variables of the
circuit.  It is shown that the Birkhoffian system associated to an
LC circuit is conservative. Under certain assumptions  on the
voltage-current characteristic for resistors, it is shown that a
Birkhoffian system associated to an RLC circuit is dissipative.
For LC/RLC networks which contain a number of loops formed only
from capacitors, the Birkhoffian associated is never regular. A
procedure to reduce the original configuration space
 to a lower dimensional one, thereby regularizing
the Birkhoffian, is  presented as well.


For RLC electrical networks, Brayton and Moser \cite{moser} proved
under a special hypothesis, that there exists a mixed potential
function which can be used to put  the system of differential
equations describing the dynamics of such a network, into a
special form (see \S 4 in \cite{moser}). The hypothesis they made
is that the currents through the inductors and the voltages across
the capacitors determine all currents and voltages in the circuit
via Kirchhoff's law.
 The
mixed potential function is constructed explicitly  only for the
networks whose graph possesses a tree containing all the capacitor
branches and none of the inductive branches, that is, the network
does not contain any loops of capacitors or cutsets of inductors,
each resistor tree branch corresponds to a current-controlled
resistor, each resistor co-tree branch corresponds to a
voltage-controlled resistor (see \S13 in \cite{moser}). Making
different assumptions on the type of admissible nonlinearities in
the circuit, this mixed potential function is used in \cite{moser}
to construct Liapunov-type
functions to prove stability.\\
 Smale \cite{smale} also develops  the differential equations for
nonlinear RLC electrical circuits and illustrates these equations
through a series of examples. He builds on the work of Brayton and
Moser \cite{moser} but he is able to treat more general equations.
A large part of the paper illustrates these equations by means of
examples and discusses stability properties of the examples.

The paper at hand is organized as follows: at the beginning of
Section 2 we present the basics of Birkhoffian systems, from the
viewpoint of differential geometry using the formalism of jets.
Then, we introduce the concepts and the direct theorems of
stability in the sense of Liapunov, within the framework of
Birkhoffian systems. In Section 3 we consider, in turn, the linear
and the nonlinear LC networks, as well as the linear and the
nonlinear RLC networks. For each of them we construct
Liapunov-type functions to prove stability or asymptotic stability
under certain conditions. Finally, we discuss an example in
Section 4.

\section{Liapunov's direct method for Birkhoffian systems}

In order to present  the ideas in a coordinate free fashion, we
consider   the formalism of 2-jets. Let $M$ be a $m$-dimensional
differentiable connected manifold. We consider the tangent bundles
($TM$, $\pi_{M}$, $M$) and ($TTM$, $\pi_{TM}$, $TM$).\\
 The \textit{2-jet bundle} ($J^2(M)$, $\pi_{J}$, $TM$) is defined
by \be J^2(M):=\left\{ z\in TTM\ /\
T\pi_{M}(z)=\pi_{TM}(z)\right\} \label{1} \ee  where
$(T\pi_M)_v:T_vTM\to T_{\pi_M(v)}M$ is the tangent map and \be
\pi_J:=\pi_{TM}\arrowvert _{J^2(M)}=T\pi_M\arrowvert_{J^2(M)}\ee A
local system of  coordinates $(q)=(q^j)_{j=1,...,m} $ on $M$
induces natural local coordinates on $J^{2}(M)$, denoted by
$(q,\dot{q},\ddot{q})=(q^j,\,\dot{q}^j,\,\ddot{q}^j)_{j=1,...,m}$
(see for example \cite{oliva}, \cite{saunders}).

A \textit{Birkhoffian} corresponding to the configuration manifold $M$
is a smooth 1-form $\omega$ on $J^2(M)$ such that, 
for any $x\in M$, we have \be \iota_{x}^*\omega=0\label{defbir}
\ee
where $\iota_x:\beta^{-1}(x)\to J^2(M)$ is the embedding of the submanifold $\beta^{-1}(x)$ into $J^2(M)$, $\beta=\pi_M\circ \pi_J$. 
From this definition it follows that, in the  natural local coordinate system ($q,\, \dot{q},\, \ddot{q}$) of $J^2(M)$, a 
Birkhoffian $\omega$ is given by \be
\omega=\sum^{m}_{j=1}\mathfrak{Q}_j(q,\, \dot{q},\,
\ddot{q})dq^j\label{bircor} \ee with certain functions
$\mathfrak{Q}_j:J^2(M)\to \mathbf{R}$. \noindent The pair ($M,\,
\omega$) is said to be a \textit{Birkhoffian system} (see
\cite{oliva}). \\
The \textit{differential system associated to a Birkhoffian}
$\omega $  is  the set (maybe empty) $D(\omega$), given by \be
D(\omega):=\left\{z\in J^2(M)\ \arrowvert\, \omega(z)=0
\right\}\label{difsistem'} \ee The manifold $M$ is the
\textit{space of configurations} of $D(\omega)$, and $D(\omega)$
is said to have $m$ 'degrees of freedom'. The $\mathfrak{Q}_i$ are
the 'generalized external forces' associated to the local
coordinate system. In the natural local coordinate system,
$D(\omega)$ is characterized by the following implicit system of
second order ODE's \be \mathfrak{Q}_j(q,\, \dot{q},\, \ddot{q})=0
\textrm{ for all } j=\overline{1,m}\label{difsistem} \ee  The
Birkhoffian formalism is  a global formalism for the dynamics of
implicit systems of second order  differential equations on a
manifold.

A cross section $X$ of the affine bundle ($J^2(M),\, \pi_J,\, TM$), that is, a smooth function $X:TM\to J^2(M)$ such that  
$\pi_J\circ X$=id,
 can be identified with a special vector field on $TM$, namely, the second
 order vector field $Y$ on $TM$,
that is, a smooth function $Y:TM\to TTM$ such that $\pi_{TM}\circ
Y$=id and  $T{\pi_M}\circ Y=$id.  Using 
the canonical   embedding  $i:J^2(M)\to TTM$,  we write  $Y=i\circ
X$.\\
 In  natural  local coordinates, a second order vector field
can be represented as \be
Y=\sum^{m}_{j=1}\left[\dot{q}^j\f{\pa}{\pa q^j}+\ddot{q}^j(q,\,
\dot{q})\f{\pa}{\pa \dot{q}^j}\right]\label{vectorfield} \ee
A \textit{Birkhoffian vector field} associated to a Birkhoffian $\omega$ of $M$ (see \cite{oliva}) is a smooth second 
order vector field on $TM$, $Y=i\circ X$, with $X:TM\to J^2(M)$,
such that $ImX\subset\, D(\omega)$, that is, \be X^*\omega=0 \ee
 In the natural  local coordinate system, a Birkhoffian vector field is given by the expression (\ref{vectorfield}), such 
that $\mathfrak{Q}_j(q,\dot{q},\ddot{q}(q,\dot{q}))=0$.

A Birkhoffian $\omega
$ is \textit{regular} if and only if \be \textrm{det}\left[\f{\pa
\mathfrak{Q}_j}{\pa \ddot{q}^i}(q,\, \dot{q},\,
\ddot{q})\right]_{j,i=1,...,m}\neq 0\label{regular} \ee for all
$(q,\, \dot{q},\, \ddot{q})$, and for each $(q,\, \dot{q})$,
there exists $(q,\dot{q},\ddot{q})\in J^2(M)$ such that $\mathfrak{Q}_j(q,\, \dot{q},\, 
\ddot{q})=0,\, j=1,...,m.$

 If  $\omega$ is  a regular Birkhoffian
corresponding to the configuration manifold $M$, then,
\textit{the principle of determinism} is satisfied, that is, there
exists an \textit{unique} Birkhoffian vector field $Y=i\circ X$
associated to $\omega$ such that $Im\, X=D(\omega)$ (see
\cite{oliva}).

A Birkhoffian $\omega$ of $M$ is called \textit{conservative} (see
\cite{oliva}) if and only if there exists a smooth function 
$E_{\omega}:TM\to \mathbf{R}$ such that \be
 (X^*\omega)Y=dE_{\omega}(Y)\label{conserv'}
\ee for all second order vector fields $Y=i\circ X$, which is
equivalent, in the natural local coordinate system, to the
identity
 \be
\sum^{m}_{j=1}\mathfrak{Q}_j(q\, \, \dot{q},\, \ddot{q})\dot{q}^j=\sum^{m}_{j=1}\left[\f{\pa E_{\omega}}{\pa q^j}\dot{q}^j+\f{\pa E_{\omega}}{\pa 
\dot{q}^j}\ddot{q}^j\right]\label{conserv} \ee If $\omega$ is
conservative and $Y$ is a Birkhoffian vector field, then
(\ref{conserv'}) becomes \be dE_{\omega}(Y)=0 \label{unic}\ee This
means that $E_{\omega}$ is constant along the trajectories of $Y$.

 A Birkhoffian $\omega$ of the configuration space $M$ is called
\textit{dissipative} (see \cite{eu2}) if and only if there exists
a smooth function $E_{0_{\omega}}:TM\to \mathbf{R}$ such that \be
(X^*\omega)Y = dE_{0_\omega}(Y)+ D(Y) \label{store1} \ee for all
second order vector fields $Y=i\circ X$ on $TM$,  $D$ being a
dissipative 1-form on $TM$, that is,
$D=\sum^{m}_{j=1}D_{j}(q,\dot{q})dq^j$ and
\be\sum^{m}_{j=1}D_{j}(q,\dot{q})\dot{q}^j>0\label{store}\ee

\noindent Equation (\ref{store1}) is  equivalent, in a local
coordinate system, to the identity \be
\sum^{m}_{j=1}\mathfrak{Q}_j(q\, \, \dot{q},\,
\ddot{q})\dot{q}^j=\sum^{m}_{j=1}\left[\f{\pa E_{0_\omega}}{\pa
q^j}\dot{q}^j+\f{\pa E_{0_\omega}}{\pa \dot{q}^j}\ddot{q}^j+
D_{j}(q,\dot{q})\dot{q}^j \right]\label{dissip} \ee

\noindent In view of (\ref{store}), we obtain from (\ref{store1}),
\be
 (X^*\omega)Y> dE_{0_{\omega}}(Y)\label{storage}
\ee for all second order vector fields $Y=i\circ X$. That is
equivalent, in local coordinates, to the dissipation inequality
\be \sum^{m}_{j=1}\mathfrak{Q}_j(q\, \, \dot{q},\,
\ddot{q})\dot{q}^j> \sum^{m}_{j=1}\left[
 \f{\pa E_{0{\omega}}}{\pa q^j}\dot{q}^j+\f{\pa E_{0_{\omega}}}{\pa
 \dot{q}^j}\ddot{q}^j\right]
\label{storage'} \ee
 If $\omega$ is a dissipative Birkhoffian  and $Y$ is the
Birkhoffian vector field, then (\ref{storage}) becomes \be
dE_{0_\omega}(Y)< 0\label{nonincreasing}\ee This means that
$E_{0_\omega}$ is nonincreasing along the trajectories of $Y$.

Let us introduce now the concepts of stability for a Birkhoffian
system.

\noindent The \textbf{\textit{equilibrium points}} of the system,
that is, the points in which the system can remain permanently at
rest, are to be found as the solutions of the system \be
\mathfrak{Q}_j(q,0,0)=0, \quad j=1,...,m\label{2} \ee Let us
denote an equilibrium point by $(q_e,0)\in\Omega\subset TM$, and
an initial state of the system by $(q_0,\dot{q}_0)$, with
$q(0)=q_0$,
$\dot{q}(0)=\dot{q}_0$. \\
For regular Birkhoffians, we can define the equilibrium points
using the notion of Birkhoffian vector field, that is, a point
$(q_e,0)$ is an equilibrium point of the Birkhoffian vector field
$Y$ if and only if \be Y(q_e,0)=0 \ee

 An equilibrium point $(q_e,0)$ is
said to be \textbf{stable} (or Liapunov stable) if for every open
neighborhood $\Omega$ of $(q_e,0)$, there is a neighborhood
$\Omega_1\subset \Omega$ such that a motion $(q(t),\dot{q}(t))$
starting  at $(q_0,\dot{q}_0)\in \Omega_1$,
 remains in $\Omega$. If in addition,
$\Omega_1$ can be chosen such that, for any $(q_0,\dot{q}_0)\in
\Omega_1$,  $(q(t),\dot{q}(t))$ converges to $(q_e,0)$
 as $t\to \infty$, then $(q_e,0)$ is said
 to be \textbf{asymptotically stable}.

In the memoir \cite{liapunov}, Liapunov presents geometric
theorems, generally referred to as the direct method of Liapunov
(see, for example, \cite{lasalle}), for deciding the stability or
instability of an equilibrium point of a differential equation.
\\
In what follows we consider \textit{\textbf{Liapunov's direct
method for Birkhoffian systems}}. This is based on finding
 a function  $V\in C^1(TM,\mathbf{R})$
such that \ba
\begin{array}{l}
\,\,(i)\,\,V(q_e,0)=0\\
\,(ii)\,V(q,\dot{q})>0 \, \textrm{ for } (q,\dot{q})\neq(q_e,0) \textrm{ in } \Omega\\
(iii)\, dV(Y)\leq 0 \, \textrm{ for all second order vector fields
$Y$ defined on
 } \Omega
 \end{array}\label{liapunovf}
\ea with $\Omega$ an open neighborhood  of $(q_e,0)$.
The function $V$ is called \textit{\textbf{Liapunov function}}.

One can prove the following theorems (completely analogous to the
theorems proved in \cite{lasalle} for a Liapunov function defined
on $U\subset M$):

\vspace{0.15cm}

\textsc{Stability theorem.} If there exists in a neighborhood
$\Omega$ of $(q_e,0)$ a Liapunov function $V(q,\dot{q})$, then
$(q_e,0)$ is stable.

\vspace{0.15cm}

\textsc{Asymptotic stability theorem.} If there exists in a
neighborhood $\Omega$ of $(q_e,0)$ a Liapunov function
$V(q,\dot{q})$, such that $dV(Y)<0$ for all second order vector
fields $Y$ defined on $\Omega$, then $(q_e,0)$ is asymptotically
stable.

\vspace{0.15cm}

From the condition ($ii$) in (\ref{liapunovf}) we get that there
exists $c_0>0$ such that the  level curve \be \{(q,\dot{q})\in
\Omega, V(q,\dot{q})=c\} \ee is a closed curve for every constant
$0\leq c\leq c_0.$ Sketching in the $m$-plane $(q,\dot{q})$ these
level curves of the function $V$,
  we obtain surfaces like "ellipsoids" centered at the equilibrium point.\\
  If $dV=0$, then the equilibrium point $(q_e,0)$ is a \textbf{center}
and the motion of the system is periodic. \\
If $dV<0$, then  each trajectory keeps moving to lower $c$ and
hence penetrates smaller and smaller "ellipsoids" as $t \to
\infty.$ Thus, the equilibrium point is asymptotically stable.
This exclude the existence of periodic motions of the system.

\section{Stability of the equilibrium points of LC and RLC
networks}

A simple electrical circuit provides us with an \textit{oriented
connected} graph. The graph will be assumed to be \textit{planar}.
Let $b$  be the total number of branches in the graph, $n$ be one
less than the number of  nodes and $m$
be the cardinality of a 
selection of loops that cover the whole graph.  By Euler's
polyhedron formula, $b=m+n$. We choose a reference node and a
current direction in each $l$-branch of the graph,
 $l=1,...,b$. We also consider a covering of the graph with $m$ loops,
   and a current direction in  each $j$-loop, 
$j=1,...,m$. We assume that the associated graph has at least one
loop, meaning that $m>0$.
 An oriented connected graph  can be described by matrices which contain
 only 0, $\pm$1, these are: the  \textit{incidence matrix}
 $B\in \mathfrak{M}_{bn}(\mathbf{R})$, rank$(B)=n$,
 and the  \textit{loop matrix} $A\in \mathfrak{M}_
{bm} (\mathbf{R})$, rank$(A)=m$. For the fundamentals of
electrical circuit theory, see, for example, \cite{chua}.

Let us now consider an RLC electrical circuit consisting of $r$
resistors, $k$ inductors and $p$ capacitors, such that to each
branch of the  associated graph there corresponds just one
electrical device, that is, $b=r+k+p$. For LC electrical circuits
$r=0$. Using the matrices $A$ and $B$, Kirchhoff's current law and
Kirchhoff's voltage law can be expressed by the equations \be B^T
\textsc{i} =0 \quad (KCL), \quad A^T v =0\quad (KVL)\label{kcl}\ee
where $\textsc{i}=(\textsc{i}_{[\Gamma]},
\textsc{i}_{(a)},\textsc{i}_{\alpha})\in \mathbf{R}^r\times
\mathbf{R}^k\times \mathbf{R}^p\simeq\mathbf{R}^b$
is the current vector
and $v=(v_{[\Gamma]},v_{(a)},v_{\alpha})\in \mathbf{R}^r\times
\mathbf{R}^k\times\mathbf{R}^p\simeq\mathbf{R}^b$
is the voltage drop vector.
 \noindent Tellegen's theorem establishes
a relation between the matrices $A^T$ and $B^T$: \textit{the
kernel of the matrix $B^T$ is orthogonal to the kernel of the
matrix $A^T$}(see, for example,  page 5 of \cite{moser}).

We consider the voltage-current laws for nonlinear devices given
by \be v_{[\Gamma]}=R_{\Gamma}(\textsc{i}_{[\Gamma]}), \quad
v_{(a)}=L_{a}(\textsc{i}_{(a)})\f{d\textsc{i}_{(a)}}{dt}, \quad
v_\alpha=C_\alpha(\textsc{q}_{\alpha}), \label{nonliniar}\ee
   $R_{\Gamma}, L_a,C_\alpha:\mathbf{R}\longrightarrow \mathbf{R}\backslash \{0\}$
being smooth functions, $\textsc{q}_\alpha$ denote the charges of
the capacitors, with
$\textsc{i}_\alpha=\f{d\textsc{q}_{\alpha}}{dt}$. If the
capacitors and the inductors are linear then the relations above
become, respectively, \be
v_{[\Gamma]}=\textrm{\scriptsize{R}}_{\Gamma}\textsc{i}_{[\Gamma]},\quad
v_{(a)}=\textrm {\scriptsize L}_a\f{d{\textsc{i}_{(a)}}}{dt}\quad
v_\alpha=\frac{\textsc{q}_\alpha}{\textrm {\scriptsize C}_\alpha},
\label{15} \ee where {\scriptsize{R}}$_{\Gamma}\neq 0$,
{\scriptsize C}$_\alpha\neq 0$ and {\scriptsize L}$_a\neq 0$ are
distinct constants.

\noindent Summing up, the equations governing the network are \be
B^T \left(
\begin{array}{c}
\textsc{i}_{[\Gamma]}\\
\textsc{i}_{(a)} \\
\f{d\textsc{q}_{\alpha}}{dt}
\end{array}
\right)=0, \quad\quad
A^T \left(
\begin{array}{c}
R_{\Gamma}(\textsc{i}_{[\Gamma]})\\
L_{a}(\textsc{i}_{(a)})\,
\f{d\textsc{i}_{(a)}}{dt} \\
C_\alpha(\textsc{q}_{\alpha})
\end{array}
\right)=0.
\label{7} \ee Using the first set of equations (\ref{7}),
 one defines (see  \cite{eu2}, \cite{io}) a  family of
  $m$-dimensional affine-linear configuration spaces
  $M_{{\bf \mathfrak{c}}}\subset\mathbf{R}^b$, parameterized by a constant vector
  $\mathfrak{c}$ in $\mathbf{R}^n$ which corresponds to initial values of certain state
variables of the circuit. Since the matrix $B$ is constant,
integrating the first set of equations (\ref{7}), one gets
$B^Tx=\mathfrak{c},$ with $\textsc{i}=\dot{x}$,
 $\mathfrak{c}$ a constant vector in $\mathbf{R}^n$. Thus, one
 defines
 \be M_{\mathfrak{c}}:=\{x\in
 \mathbf{R}^b|B^Tx=\mathfrak{c}\}\label{20'}
\ee  Its dimension is $m=b-n$, because rank$(B)=n$. Local
coordinates on $M_\mathfrak{c}$ are denoted by $q=(q^1,..,q^m)$.
Solving the system in (\ref{20'}),
 one expresses any of the $x$-variables
in terms  of $q$-s, namely, \be x=\mathcal{N}q+\mathcal{K}
\label{n}\ee where $\mathcal{N}=\left(\begin{array}{c}
\mathcal{N}^{\Gamma}_j\\
\mathcal{N}^{a}_{ j}\\
\mathcal{N}^{\alpha}_{j}
\end{array}\right)_{\Gamma=\overline{1,r},a=\overline{r+1,r+k},
\alpha=\overline{r+k+1,b},\atop j=\overline{1,m}}$ is a matrix of
constants and $\mathcal{K}=\left(\begin{array}{c}
\mathcal{K}^{\Gamma}\\
\mathcal{K}^{a}\\
\mathcal{K}^{\alpha}
\end{array}\right)
$
a constant
vector in $\mathbf{R}^b$.\\
  By Tellegen's theorem and a
fundamental theorem of linear algebra, one obtains that
$Ker(A^T)=Ker(\mathcal{N}^T)$ (see \cite{eu2} , \cite{io}). \\
A
Birkhoffian $\omega_{\mathfrak{c}}$ on the configuration space
$M_{\mathfrak{c}}$ arises from a linear combination of the second
set of equations (\ref{7}), by replacing the matrix $A^T$ with the
matrix of constants $\mathcal{N}^T$.

\vspace{0.3cm}

 I) For a linear LC
network ($r=0$) we have the following expression of the
Birkhoffian (see \cite{io}) \be
\mathfrak{Q}_j(q,\dot{q},\ddot{q})=\sum_{a=1}^{k}\sum_{i=1}^{m}
{\textrm {\scriptsize L}}_a
\mathcal{N}^{a}_j\mathcal{N}^{a}_i\ddot{q}^i+ \sum
^{b}_{\alpha=k+1} \sum^{m}_{i=1}\f{\mathcal{N}_j^{\,
\alpha}\mathcal{N}^ {\alpha}_{\,  i}}{{\textrm {\scriptsize
C}}_{\alpha-k}}q^i+(\textrm{const})_j \label{4} \ee with const
$\in\mathbf{R}^m$ a constant vector. \\
A linear LC network is \textit{conservative} (see \cite{io}). The
function $E_{\omega}:TM_\mathfrak{c}\rightarrow \mathbf{R}$
satisfying (\ref{conserv}) has the following expression
 \ba E_{\omega}(q,\dot{q})&=&\f{1}{2}\sum^{k}_{a=1}\sum_{j,i=1}^m
{\textrm {\scriptsize L}}_a
\mathcal{N}^{a}_j\mathcal{N}^{a}_i\dot{q}^j \dot{q}^i+\f{1}{2}
\sum^{b}_{\alpha=k+1}\sum^{m}_{j,i=1}
\f{\mathcal{N}^{\alpha}_j\mathcal{N}^{\alpha}_i} {{\textrm
{\scriptsize
C}}_{\alpha-k}}q^jq^i\nonumber\\
&& + \sum_{j=1}^{m}(\textrm{const})_jq^j \label{energie} \ea

 In what follows we assume that \be
\textrm{det}\left[\sum^{k}_{a=1}{\textrm {\scriptsize
L}}_a\mathcal{N}_j^{a}\mathcal{N}^{a}_{i}\right]_{j,i=1,...,m}
\neq 0, \quad
\textrm{det}\left[\sum^{b}_{\alpha=k+1}\f{\mathcal{N}_j^{\,
\alpha}\mathcal{N}^ {\alpha}_{\,  i}}{{\textrm {\scriptsize
C}}_{\alpha-k}}\right]_{j,i=1,...,m}\neq 0\label{det}\ee
 that is, the network does not contain loops formed only by
 capacitors and respectively, loops formed only by
 inductors (see \cite{io}). If the network contains
  capacitor
loops and  inductor loops,
 we will reduce  first the configuration space to a lower dimensional configuration
 space. On the
reduced configuration space the corresponding Birkhoffian is still
conservative (see \cite{io}) and  the corresponding determinants
(\ref{det}) will be different from zero.
 The inductor loops can be
considered as some conserved quantities of the network.

\vspace{0.15cm}

\textsc{Theorem 1.} \textit{Let $(q_e,0)$ be an equilibrium point
of a linear LC network with the Birkhoffian components given by
(\ref{4}). Then $q_e$ satisfies the system }\be \sum
^{b}_{\alpha=k+1} \sum^{m}_{i=1}\f{\mathcal{N}_j^{\,
\alpha}\mathcal{N}^ {\alpha}_{\,  i}}{{\textrm {\scriptsize
C}}_{\alpha-k}}q^i+(\textrm{const})_j=0, \quad j=1,...,m\label{5}
\ee \textit{To each} $\textrm{const}$ \textit{which is related to
the initial data for the considered network,
 we get \textbf{a unique
equilibrium point}. If}
 \be
 {\textrm {\scriptsize L}}_a>0,\,
\forall \, a=1,...,k, \quad {\textrm {\scriptsize
C}}_{\alpha}>0,\, \forall \, \alpha= 1,...,p\label{condliniar}\ee
\textit{the equilibrium point is a \textbf{stable center}, and the
motion of the system is periodic.}

\vspace{0.15cm}

Indeed, the equilibrium points of a linear LC network are obtained
as solutions of the system $\mathfrak{Q}_j(q,0,0)=0$, $j=1,...,m$,
where $\mathfrak{Q}_j(q,\dot{q},\ddot{q})$ is given by (\ref{4}).
Thus, we get that $q_e$ has to fulfill the system (\ref{5}). Under
the second condition in (\ref{det}), this system has
 for each const $\in \mathbf{R}^m$
a unique solution.\\
 The stability of this equilibrium point is
obtained by the Stability Theorem presented in section 2. We
define a Liapunov function $V\in C^1(TM_\mathfrak{c},\mathbf{R})$
by \ba \hspace{-0.8cm}
V(q,\dot{q})&=&E_\omega(q,\dot{q})-E_\omega(q_e,0)\label{V}
\nonumber\\
&=&\f{1}{2}\sum^{k}_{a=1}\sum_{j,i=1}^m {\textrm {\scriptsize
L}}_a \mathcal{N}^{a}_j\mathcal{N}^{a}_i\dot{q}^j
\dot{q}^i+\f{1}{2} \sum^{b}_{\alpha=k+1}\sum^{m}_{j,i=1}
\f{\mathcal{N}^{\alpha}_j\mathcal{N}^{\alpha}_i} {{\textrm
{\scriptsize C}}_{\alpha-k}}(q^j-q^j_e)(q^i-q_e^i) \nonumber\\
&&\ea where $q_e$ satisfies the system (\ref{5}). Indeed, this
function satisfies the conditions (\ref{liapunovf}). Taking into
account (\ref{condliniar}), the  matrices $\left(\sum ^{k}_{a=1}
{\textrm {\scriptsize L}}_a \mathcal{N}^{a}_j\mathcal{N}^{a}_i
\right)_{j,i}$ and $\left(\sum ^{b}_{\alpha=k+1}
\f{\mathcal{N}_j^{\, \alpha}\mathcal{N}^ {\alpha}_{\, i}}{{\textrm
{\scriptsize C}}_{\alpha-k}}\right)_{j,i}$
  are positive definite. Thus, the condition $(ii)$  in
  (\ref{liapunovf}) is fulfilled. The first determinant in
(\ref{det}) being different from zero implies that the
corresponding Birkhoffian is regular. Therefore, along the
trajectories of the unique (principle of determinism) Birkhoffian
vector field, the function $E_{\omega}$  defined in
(\ref{energie}), satisfies (\ref{unic}). Thus, the function
(\ref{V}) satisfies the condition $(iii)$ in (\ref{liapunovf}).
 In this case, sketching in the  $m$-plane
$(q,\dot{q})$ the level curves of the function (\ref{V}),
  we obtain  ellipsoids centered at the equilibrium point.
  The equilibrium point is a center and the motion of the system is
  periodic.$\square$
\\

II) For a nonlinear LC network we have the following expression of
the Birkhoffian (see \cite{io})  \ba
\hspace{-0.6cm}\mathfrak{Q}_j(q, \, \dot{q},\,\ddot{q})&=
&\sum^{k}_{a=1}\mathcal{N}^{\,  a}_{j} L_a\left(
\sum^{m}_{l=1}\mathcal{N}^ {a}_{\, l}\dot{q}^l\right)\sum
^{m}_{i=1}\mathcal{N}^ {a}_{\, i}\ddot{q}^i+ \sum
^{b}_{\alpha=k+1}\mathcal{N}_j^{\,
\alpha}C_{\alpha-k}\left(\sum^{m}_{l=1}\mathcal{N}^ {\alpha}_{\,
l}q^l+\mathcal{K}^\alpha\right) \nonumber\\
&=&\sum^{m}_{i=1} \left(\sum^{k}_{a=1}\mathcal{N}_j^{\,
a}\mathcal{N}^ {a}_{\,
i}\widetilde{L}_a\left(\dot{q}\right)\right)\ddot{q}^i + \sum
^{b}_{\alpha=k+1}\mathcal{N}_j^{\,
\alpha}\widetilde{C}_{\alpha-k}\left(q\right) \label{bir} \ea

\noindent A nonlinear LC network is \textit{conservative} (see
\cite{io}). In this case, the function
$E_{\omega}:TM_\mathfrak{c}\rightarrow \mathbf{R}$ is given by \be
E_{\omega}(q,\dot{q})=\mathcal{E}(\dot{q})+\mathfrak{E}(q)
\label{energienonlinear} \ee with \ba \hspace{-0.6cm}
\mathcal{E}(\dot{q})&=&\sum_{a=1}^{k}\sum_{l=1}^m\sum_{i_1<...<i_l=1}^{m}(-1)^{l+1}\underbrace{\int_{...}
\int}_l \left[ \widetilde{L}_a^{(l-1)}(\dot{q})\mathcal{N}^{
a}_{i} \dot{q}^i\right.\nonumber\\
\hspace{-0.6cm}&& \hspace{0.8cm}\left. \,\, \quad \quad \quad
\quad \quad  \quad \quad \quad +\,
(l-1)\widetilde{L}_a^{(l-2)}(\dot{q})\right]\mathcal{N}^{
a}_{i_1}...\mathcal{N}^{ a}_{i_l}
d\dot{q}^{i_1}...d\dot{q}^{i_l}\nonumber\\
\hspace{-0.6cm}\mathfrak{E}(q)&=&\sum_{\alpha=k+1}^{b}\sum_{l=1}^m\sum_{i_1<...<i_l=1}^{m}(-1)^{l+1}\underbrace{\int_{...}
\int}_l
\widetilde{C}_{\alpha-k}^{(l-1)}(q)\mathcal{N}^{\alpha}_{i_1}...\mathcal{N}^{\alpha}_{i_l}
dq^{i_1}...dq^{i_l} \label{6} \ea
where we denoted the derivatives $\widetilde{C}^{(l)}_{\alpha-k}:=\f {d^l\widetilde{C}_{\alpha}(x)}{dx^l}$, $\widetilde{L}^{(l)}_{a}:=\f 
{d^l\widetilde{L}_a(x)}{dx^l}$.

 In what follows we assume that
\be \textrm{det}\left[\sum^{k}_{a=1}
\mathcal{N}_j^{a}\mathcal{N}^{a}_{i}\widetilde{L}_a(\dot{q})\right]_{j,i=1,...,m}
\neq 0\label{det2}\ee
 that is, the network does not contain capacitor loops.
 In the case the network contains capacitor loops,
 we first reduce  the configuration
space to a lower dimensional one. On the reduced configuration
space the corresponding Birkhoffian is still conservative (see
\cite{io}) and  the corresponding determinant above will be
different from zero.

\vspace{0.15cm}

\textsc{Theorem 2.} \textit{Let $(q_e,0)$ be an equilibrium point
of a nonlinear LC network with the Birkhoffian components given by
(\ref{bir}). Then $q_e$ satisfies the system } \be \sum
^{b}_{\alpha=k+1}\mathcal{N}_j^{\,
\alpha}C_{\alpha-k}\left(\sum^{m}_{l=1}\mathcal{N}^ {\alpha}_{\,
l}q^l+\mathcal{K}^\alpha\right)=0, \quad j=1,...,m\label{5''} \ee
\textit{A nonlinear LC networks can have \textbf{several
equilibrium points}. If \be L_a(0)>0,\,  \forall \, a=1,...,k,
\quad C'_{\alpha}(q_e)>0,\,
 \forall \, \alpha=1,...,p \label{cond}\ee then  the equilibrium points are
\textbf{locally stable centers}.}

\vspace{0.15cm}

Indeed, the equilibrium points of a nonlinear LC network are
obtained as solutions of the system $\mathfrak{Q}_j(q,0,0)=0$,
$j=1,...,m$, where $\mathfrak{Q}_j(q,\dot{q},\ddot{q})$ is given
by (\ref{bir}). Thus, we see that $q_e$ has to fulfill the
system (\ref{5''}).\\
 The local stability of the equilibrium points is
obtained by the Stability Theorem presented in section 2. We
define a Liapunov function $V\in C^1(TM_\mathfrak{c},\mathbf{R})$
by \be
V(q,\dot{q})=E_\omega(q,\dot{q})-E_\omega(q_e,0)\label{V'}\ee with
$E_\omega$ given by (\ref{energienonlinear}) and $q_e$ satisfying
the system (\ref{5''}). \\
 Let us now evaluate the Hessian matrix
of the function $V$ in (\ref{V'}) at an equilibrium point
$(q_e,0)$. We get \be \textbf{H}_{V}(q_e,0)=
\left(\begin{array}{cc} \f{\pa^2 \mathcal{E}(\dot{q})}{\pa
\dot{q}^i\pa \dot{q}^j}|_{(q_e,0)} & 0
\\
0 & \f{\pa^2 \mathfrak{E}(q)}{\pa q^i\pa q^j}|_{(q_e,0)}
\end{array}\right)\label{hes}
\ee For the Birkhoffian (\ref{bir}), the function $E_{\omega}$ in
(\ref{energienonlinear}) satisfies the identity (\ref{conserv})
(see \cite{io}), that is, \ba \f{\pa \mathcal{E}(\dot{q})}{\pa
\dot{q}^i}&=&\sum_{a=1}^{k}\sum_{l=1}^m\widetilde{L}_a(\dot{q})\mathcal{N}^{
a}_{i} \mathcal{N}^{ a}_{l}\dot{q}^l
\\
\f{\pa^2 \mathfrak{E}(q)}{\pa q^i}&=&\sum_{\alpha=k+1}^{b}
\widetilde{C}_{\alpha-k}(q)\mathcal{N}^{\alpha}_{i} \ea Therefore,
we get
 \ba \f{\pa^2 \mathcal{E}(\dot{q})}{\pa
\dot{q}^i\pa \dot{q}^j}&=&\sum_{a=1}^{k}\left[
\sum_{l=1}^m\widetilde{L}_a'(\dot{q})\mathcal{N}^{ a}_{j}
\mathcal{N}^{ a}_{i}\mathcal{N}^{ a}_{l}\dot{q}^l+
\widetilde{L}_a(\dot{q})\mathcal{N}^{ a}_{i}\mathcal{N}^{ a}_{j}
\right]\label{mathcalE}
\\
\f{\pa^2 \mathfrak{E}(q)}{\pa q^i\pa q^j}&=&\sum_{\alpha=k+1}^{b}
\widetilde{C}'_{\alpha-k}(q)\mathcal{N}^{\alpha}_{i}
\mathcal{N}^{\alpha}_{j}\label{mathfrakE} \ea From
(\ref{mathcalE}), (\ref{mathfrakE}), the matrix (\ref{hes}) writes
as \be \textbf{H}_{V}(q_e,0)= \left(\begin{array}{cc}
\sum_{a=1}^{k}L_a(0)\mathcal{N}^{ a}_{i}\mathcal{N}^{ a}_{j}& 0\\
0 & \sum_{\alpha=k+1}^{b}
\widetilde{C}'_{\alpha-k}(q_e)\mathcal{N}^{\alpha}_{i}\mathcal{N}^{\alpha}_{j}
\end{array}\right)\label{hessiana}
\ee

\noindent In view of the conditions (\ref{cond}), the matrices
$\left( \sum_{a=1}^{k}L_a(0)\mathcal{N}^{
a}_{i}\mathcal{N}^{ a}_{j} \right)_{i,j}$ and\\
$\left(\sum_{\alpha=k+1}^{b}
\widetilde{C}'_{\alpha-k}(q_e)\mathcal{N}^{\alpha}_{i}\mathcal{N}^{\alpha}_{j}\right)_{i,j}$
 are  positive definite. Therefore,  the Hessian matrix
 (\ref{hessiana})
is positive definite. The centers of the level curves of the
function (\ref{V'}) have the coordinates $(q_e,0)$, where $q_e$
satisfies the system (\ref{5''}). Thus, in a neighborhood of an
equilibrium point, the condition $(ii)$ in
  (\ref{liapunovf}) is fulfilled by the function $V$ in (\ref{V'}).
   The determinant (\ref{det2}) being different
from zero implies that the corresponding Birkhoffian is regular.
Therefore, along the trajectories of the unique (principle of
determinism) Birkhoffian vector field, the function $E_{\omega}$
satisfies (\ref{unic}). Thus, the function (\ref{V'}) satisfies
the condition $(iii)$ in (\ref{liapunovf}).
  By the Stability Theorem, the equilibrium points are locally
  stable centers.$\square$
\\

III) For a linear RLC network we have the following expression of
the Birkhoffian (see \cite{eu2}) \ba
\mathfrak{Q}_j(q,\dot{q},\ddot{q})&=&\sum_{a=r+1}^{r+k}\sum_{i=1}^{m}
{\textrm {\scriptsize L}}_{a-r}
\mathcal{N}^{a}_j\mathcal{N}^{a}_i\ddot{q}^i+ \sum
^{r}_{\Gamma=1}\sum^{m}_{i=1}\textrm {\scriptsize
R}_{\Gamma}\mathcal{N}_{j}^{\Gamma} \mathcal{N}^
{\Gamma}_{i}\dot{q}^i\nonumber\\
&&+ \sum ^b_{\alpha=r+k+1} \sum^{m}_{i=1}\f{\mathcal{N}_j^{\,
\alpha}\mathcal{N}^ {\alpha}_{\,  i}}{{\textrm {\scriptsize
C}}_{\alpha-r-k}}q^i+(\textrm{const})_j  \label{4'} \ea with const
$\in\mathbf{R}^m$ a constant vector. \\
A linear RLC network with \be \textrm {\scriptsize
R}_{\Gamma}>0,\quad
 \Gamma=1,...,r,
 \label{R}\ee is \textit{dissipative} (see \cite{eu2}). The
function $E_{0_\omega}:TM_\mathfrak{c}\rightarrow \mathbf{R}$ and
the dissipative 1-form  satisfying (\ref{dissip}), are given by
\ba
E_{0_\omega}(q,\dot{q})&=&\f{1}{2}\sum^{r+k}_{a=r+1}\sum_{j,i=1}^m
{\textrm {\scriptsize L}}_{a-r}
\mathcal{N}^{a}_j\mathcal{N}^{a}_i\dot{q}^j \dot{q}^i+\f{1}{2}
\sum^b_{\alpha=r+k+1}\sum^{m}_{j,i=1}
\f{\mathcal{N}^{\alpha}_j\mathcal{N}^{\alpha}_i} {{\textrm
{\scriptsize
C}}_{\alpha-r-k}}q^jq^i\nonumber\\
&&+\sum_{j=1}^{m}(\textrm{const})_jq^j \label{dissipative} \ea \be
 D=\sum^{m}_{j,i=1}\sum ^{r}_{\Gamma=1}\textrm {\scriptsize
R}_{\Gamma}\mathcal{N}_{j}^{\Gamma} \mathcal{N}^
{\Gamma}_{i}\dot{q}^idq^j\ee

In what follows we assume that \be
\hspace{-0.7cm}\textrm{det}\left[\sum^{r+k}_{a=r+1}{\textrm
{\scriptsize
L}}_{a-r}\mathcal{N}_j^{a}\mathcal{N}^{a}_{i}\right]_{j,i=1,...,m}
\neq 0, \quad
\textrm{det}\left[\sum^{b}_{\alpha=r+k+1}\f{\mathcal{N}_j^{\,
\alpha}\mathcal{N}^ {\alpha}_{\,  i}}{{\textrm {\scriptsize
C}}_{\alpha-r-k}}\right]_{j,i=1,...,m}\neq 0\label{det'}\ee
 that is, the network does not contain capacitor loops and
 inductor loops, respectively. If the network contains
 capacitor
loops and inductor loops,  we will reduce first the configuration
space to a lower dimensional configuration space. On the reduced
configuration space the corresponding Birkhoffian is still
dissipative (see \cite{eu2}) and the corresponding determinants
above will be different from zero.

\vspace{0.15cm}

\textsc{Theorem 3.} \textit{Let $(q_e,0$) be an equilibrium point
of a linear RLC network with the Birkhoffian given by (\ref{4'}).
Then $q_e$ satisfies the system} \be \sum ^b_{\alpha=r+k+1}
\sum^{m}_{i=1}\f{\mathcal{N}_j^{\, \alpha}\mathcal{N}^
{\alpha}_{\,  i}}{{\textrm {\scriptsize
C}}_{\alpha-r-k}}q^i+(\textrm{const})_j=0, \quad
j=1,...,m\label{5'} \ee \textit{To each}  $\textrm{const}$
\textit{which is related to the initial data for the considered
network,
 we get \textbf{a unique
equilibrium point}}. \textit{If \be {\textrm {\scriptsize
L}}_a>0,\, \forall \, a=1,...,k,\quad {\textrm {\scriptsize
C}}_{\alpha}>0,\, \forall \, \alpha= 1,...,p\label{condliniar'}\ee
the equilibrium point is \textbf{asymptotically stable}}.

\vspace{0.15cm}

\noindent Indeed, the equilibrium points of a linear RLC network
are obtained as solutions of the system
$\mathfrak{Q}_j(q,0,0)=0,\, j=1,...,m$, where
$\mathfrak{Q}_j(q,\dot{q},\ddot{q})$ is given by (\ref{4'}). Thus,
we see that $q_e$ has to fulfill  the system (\ref{5'}). Under the
second condition in (\ref{det'}), this system has for each const
$\in \mathbf{R}^m$ a unique
solution.\\
The asymptotic stability of this equilibrium point is obtained by
the Asymptotic Stability Theorem presented in section 2. We define
a Liapunov function $V\in C^1(TM_\mathfrak{c},\mathbf{R})$ by \ba
 V(q,\dot{q})&=&E_{0_\omega}(q,\dot{q})-E_{0_\omega}(q_e,0) =
\f{1}{2}\sum^{r+k}_{a=r+1}\sum_{j,i=1}^m {\textrm {\scriptsize
L}}_{a-r} \mathcal{N}^{a}_j\mathcal{N}^{a}_i\dot{q}^j
\dot{q}^i\nonumber\\
&&+\f{1}{2} \sum^{b}_{\alpha=r+k+1}\sum^{m}_{j,i=1}
\f{\mathcal{N}^{\alpha}_j\mathcal{N}^{\alpha}_i} {{\textrm
{\scriptsize C}}_{\alpha-r-k}}(q^j-q^j_e)(q^i-q_e^i)\label{V''}
\ea where $q_e$ satisfies the system (\ref{5'}). Indeed, this
function satisfies the conditions (\ref{liapunovf}). Taking into
account (\ref{condliniar'}), the  matrices $\left(\sum
^{r+k}_{a=r+1} {\textrm {\scriptsize L}}_{a-r}
\mathcal{N}^{a}_j\mathcal{N}^{a}_i \right)_{j,i}$ and $\left(\sum
^{b}_{\alpha=r+k+1} \f{\mathcal{N}_j^{\, \alpha}\mathcal{N}^
{\alpha}_{\, i}}{{\textrm {\scriptsize
C}}_{\alpha-r-k}}\right)_{j,i}$
  are positive definite. Thus, the condition $(ii)$  in
  (\ref{liapunovf}) is fulfilled. The first determinant in
(\ref{det'}) being different from zero implies that the
corresponding Birkhoffian is regular. Therefore, along the
trajectories of the unique (principle of determinism) Birkhoffian
vector field, the function $E_{0_\omega}$  satisfies
(\ref{nonincreasing}). Thus, the function (\ref{V''}) also
satisfies (\ref{nonincreasing}). In this case, sketching  the
level curves of the function (\ref{V''}) in the  $m$-plane
$(q,\dot{q})$,
  we obtain  ellipsoids centered at the equilibrium point.
From the Asymptotic Stability Theorem we conclude that the
equilibrium point is asymptotically stable. This excludes the
existence of periodic motions of the system.$\square$
\\

IV) For a nonlinear RLC network we have the following expression
of the Birkhoffian  (see \cite{eu2}) \ba
\hspace{-0.8cm}\mathfrak{Q}_j(q, \, \dot{q},\,\ddot{q})&=&
\sum^{r+k}_{a=r+1}\mathcal{N}^a_{j}L_{a-r}\left(\sum^{m}_{l=1}\mathcal{N}^
{a}_{l}\dot{q}^l\right)\left(\sum ^{m}_{i=1}\mathcal{N}^
{a}_{i}\ddot{q}^i\right)+\nonumber\\
&&\sum
^{r}_{\Gamma=1}\mathcal{N}_j^{\Gamma}R_{\Gamma}\left(\sum^{m}_{l=1}\mathcal{N}^
{\Gamma}_{l}\dot{q}^l\right)+\sum
^{b}_{\alpha=r+k+1}\mathcal{N}_j^{\alpha}C_{\alpha-r-k}\left(\sum^{m}_{l=1}\mathcal{N}^
{\alpha}_{l}q^l+\mathcal{K}^\alpha\right)\nonumber\\
 &=& \sum^{m}_{i=1}
\sum^{r+k}_{a=r+1}\mathcal{N}^a_{j}\mathcal{N}^
{a}_{i}\widetilde{L}_{a-r}\left(\dot{q}\right)\ddot{q}^i+ \sum
^{r}_{\Gamma=1}\mathcal{N}_j^{\Gamma}\widetilde{R}_{\Gamma}\left(\dot{q}\right)+
\sum
^{b}_{\alpha=r+k+1}\mathcal{N}_j^{\alpha}\widetilde{C}_{\alpha-r-k}\left(q\right)
\nonumber\\
&& \label{birRLC} \ea

 \noindent In order to obtain  a
\textit{dissipative} Birkhoffian (see \cite{eu2}), we assume that,
for all $x \neq 0$, \be xR_{\Gamma}(x)>0, \quad
 \forall\, \Gamma=1,...,r\label{assumption} \ee that is, for each
nonlinear resistor, the graph of the function $R_{\Gamma}$ lies in
the first and  third quadrants. The function
$E_{0_\omega}:TM_\mathfrak{c}\rightarrow \mathbf{R}$ and the
dissipative 1-form  satisfying (\ref{dissip}), are given by
 \be
E_{0_\omega}(q,\dot{q})=\mathcal{E}_0(\dot{q})+\mathfrak{E}_0(q)
\label{energienonlinear'} \ee with \ba
\hspace{-0.7cm}\mathcal{E}_0(\dot{q})&=&\sum_{a=r+1}^{r+k}
\sum_{l=1}^m\sum_{i_1<...<i_l=1}^{m}(-1)^{l+1}\underbrace{\int_{...}
\int}_l \left[ \widetilde{L}_{a-r}^{(l-1)}(\dot{q})\mathcal{N}^{
a}_{i} \dot{q}^i\right.\nonumber\\
\hspace{-0.7cm}&&\hspace{1cm}
 \left. \,\, \quad \quad \quad  \quad \quad  \quad \quad  \quad
+\, (l-1)\widetilde{L}_{a-r}^{(l-2)}(\dot{q})\right]\mathcal{N}^{
a}_{i_1}...\mathcal{N}^{ a}_{i_l}
d\dot{q}^{i_1}...d\dot{q}^{i_l}\nonumber\\
\hspace{-0.7cm}\mathfrak{E}_0(q)&=&\sum_{\alpha=r+k+1}^{b}\sum_{l=1}^m\sum_{i_1<...<i_l=1}^{m}(-1)^{l+1}\underbrace{\int_{...}
\int}_l
\widetilde{C}_{\alpha-r-k}^{(l-1)}(q)\mathcal{N}^{\alpha}_{i_1}...\mathcal{N}^{\alpha}_{i_l}
dq^{i_1}...dq^{i_l}\label{6'}\ea and \be D=\sum^{m}_{j=1}\sum
^{r}_{\Gamma=1}\mathcal{N}_j^{\Gamma}\widetilde{R}_{\Gamma}\left(\dot{q}\right)dq^j\ee

In what follows we assume that \be
\textrm{det}\left[\sum^{r+k}_{a=r+1}
\widetilde{L}_{a-r}(\dot{q})\mathcal{N}_j^{a}\mathcal{N}^{a}_{i}
\right]_{j,i=1,...,m} \neq 0\label{det2'}\ee
 that is, the network does not contain
  capacitor loops. In the case that the network contains capacitor loops,
   first, we first reduce the configuration space to a
lower dimensional configuration space on which the corresponding
Birkhoffian is still dissipative (see \cite{eu2}) and on which the
corresponding determinant above will be different from zero.

\vspace{0.15cm}

\textsc{Theorem 4.} \textit{Let $(q_e,0$) be an equilibrium point
of a nonlinear RLC network with the Birkhoffian given by
(\ref{birRLC}). Then $q_e$ satisfies the system} \be
\hspace{-0.7cm} \sum
^{r}_{\Gamma=1}\mathcal{N}_j^{\Gamma}R_{\Gamma}\left(0\right)
+\sum ^{b}_{\alpha=r+k+1}\mathcal{N}_j^{\,
\alpha}C_{\alpha-r-k}\left(\sum^{m}_{l=1}\mathcal{N}^ {\alpha}_{\,
l}q^l+\mathcal{K}^\alpha\right)=0, \quad j=1,...,m\label{sys2}\ee
\textit{A nonlinear RLC network can have  \textbf{several
equilibrium points}}.

\noindent 1) \textit{If \be R_{\Gamma}(0)=0,\quad \forall \,
\Gamma=1,...,r \label{r}\ee \be L_a(0)>0,\,  \forall \, a=1,...,k,
\quad C'_{\alpha}(q_e)>0,\,
 \forall \, \alpha=1,...,p \label{cond'}\ee
the equilibrium points are \textbf{locally asymptotically
stable}.}
\\

 \noindent 2) \textit{ If there exists $\Gamma=1,...,r$
such that $R_{\Gamma}(0)\neq 0$, but, for all $x\neq0$, \be
x\left(R_{\Gamma}(x)-R_{\Gamma}(0)\right)
>0, \quad \forall \, \Gamma=1,...,r \label{second}\ee
 and the conditions (\ref{cond'}) are fulfilled,
 then
the equilibrium points are \textbf{locally asymptotically
stable}.}

\vspace{0.15cm}

Indeed, the equilibrium points of a nonlinear RLC network are
obtained as solutions of the system $\mathfrak{Q}_j(q,0,0)=0$,
$j=1,...,m$, where $\mathfrak{Q}_j(q,\dot{q},\ddot{q})$ is given
by (\ref{birRLC}). Thus, we see that $q_e$ has to fulfill the
system (\ref{sys2}).
 The local asymptotic stability of the equilibrium points now follows from
 the Asymptotic Stability Theorem presented in section
2.\\
First we assume condition (\ref{r}) to be satisfied. Then the
system (\ref{sys2}) writes as \be \sum
^{b}_{\alpha=r+k+1}\mathcal{N}_j^{\,
\alpha}C_{\alpha-r-k}\left(\sum^{m}_{l=1}\mathcal{N}^ {\alpha}_{\,
l}q^l+\mathcal{K}^\alpha\right)=0, \quad j=1,...,m\label{sys2'}\ee
 In order
to show 1), we define a Liapunov function $V\in
C^1(TM_\mathfrak{c},\mathbf{R})$ by \be
V(q,\dot{q})=E_{0_\omega}(q,\dot{q})-E_{0_\omega}(q_e,0)\label{V'''}\ee
with $E_{0_\omega}$ given by (\ref{energienonlinear'}) and $q_e$
satisfying the system (\ref{sys2'}). In the neighborhood of any
equilibrium point,  this function satisfies the conditions
(\ref{liapunovf}). Taking into account the conditions
(\ref{cond'}),
 the Hessian
matrix of the function $V$ in (\ref{V'''}), at the equilibrium
point $(q_e,0)$, \be \textbf{H}_{V}(q_e,0)=
\left(\begin{array}{cc}
\sum_{a=r+1}^{r+k}L_{a-r}(0)\mathcal{N}^{ a}_{i}\mathcal{N}^{ a}_{j}& 0\\
0 & \sum_{\alpha=r+k+1}^{b}
\widetilde{C}'_{\alpha-r-k}(q_e)\mathcal{N}^{\alpha}_{i}\mathcal{N}^{\alpha}_{j}
\end{array}\right)\label{hessiana'}
\ee
 is positive
definite. The centers of the level curves of the function
(\ref{V'''}) have the coordinates $(q_e,0)$, where $q_e$ satisfies
the system (\ref{sys2'}). Thus, in a neighborhood of an
equilibrium point the condition ($ii$) in (\ref{liapunovf}) is
fulfilled by the function $V$ in (\ref{V'''}).
 The determinant (\ref{det2'}) being different from
zero implies that the corresponding Birkhoffian is regular.
Therefore, along the trajectories of the unique (principle of
determinism) Birkhoffian vector field, the function $E_{0_\omega}$
satisfies (\ref{nonincreasing}). Thus, the function (\ref{V'''})
also satisfies (\ref{nonincreasing}).  By the Asymptotic Stability
Theorem, the equilibrium points are locally asymptotically
stable.\\
 We assume now that there exists $\Gamma=1,...,r$
such that $R_{\Gamma}(0)\neq 0$. In order to show 2),  we consider
instead of  the function $E_{0_\omega}$  the following function
{\Large E}$_{0_\omega}:TM_\mathfrak{c} \rightarrow \mathbf{R}$ \be
\textrm{\Large
E}_{0_\omega}(q,\dot{q})=\mathcal{E}_0(\dot{q})+\mathfrak{E}_0(q)+\sum^r_{\Gamma=1}\mathcal{N}_j^{\Gamma}R_{\Gamma}(0)q^j
\label{energienonlinear''} \ee $\mathcal{E}_0(\dot{q})$,
$\mathfrak{E}_0(q)$ being given by (\ref{6'}), and instead of $D$
 the following dissipative 1-form \be \textrm{\Large
D}=\sum^{m}_{j=1}\sum
^{r}_{\Gamma=1}\mathcal{N}_j^{\Gamma}\left[\widetilde{R}_{\Gamma}
\left(\dot{q}\right)-R_{\Gamma}(0)\right]dq^j\label{seconddissip}\ee
\noindent In  view of assumption (\ref{second}), the vertical
1-form (\ref{seconddissip}) is indeed dissipative, that is, \be
\sum^{m}_{j=1}\sum
^{r}_{\Gamma=1}\left(\mathcal{N}_j^{\Gamma}\dot{q}^j\right)\left[R_{\Gamma}\left(\sum^{m}_{l=1}\mathcal{N}^
{\Gamma}_{l}\dot{q}^l\right)-R_{\Gamma}(0)\right]>0\ \ee One can
easily check that  for the function $\textrm{\Large
E}_{0_\omega}(q,\dot{q})$ given by (\ref{energienonlinear''}) and
the dissipative 1-form (\ref{seconddissip}), the Birkhoffian
(\ref{birRLC}) is dissipative, that is, the identity \be
\sum^{m}_{j=1}\mathfrak{Q}_j(q\, \, \dot{q},\,
\ddot{q})\dot{q}^j=\sum^{m}_{j=1}\left[\f{\pa \textrm{\Large
E}_{0_\omega}}{\pa q^j}\dot{q}^j+\f{\pa \textrm{\Large
E}_{0_\omega}}{\pa \dot{q}^j}\ddot{q}^j+ \textrm{\Large
D}_{j}(q,\dot{q})\dot{q}^j \right] \ee
is fulfilled. \\
We define now a Liapunov function $\textrm{\Large V}\in
C^1(TM_\mathfrak{c},\mathbf{R})$ by \be \textrm{\Large
V}(q,\dot{q})=\textrm{\Large
E}_{0_\omega}(q,\dot{q})-\textrm{\Large
E}_{0_\omega}(q_e,0)\label{V''''}\ee If $\omega$ is a dissipative
Birkhoffian  and $Y$ is the Birkhoffian vector field, then
(\ref{nonincreasing}) becomes \be d\textrm{\Large
E}_{0_\omega}(Y)< 0\label{nonincreasing'}\ee The function
$\textrm{\Large V}$ in (\ref{V''''}) satisfies
(\ref{nonincreasing'}) as well.\\
The centers of the level curves of the function (\ref{V''''}) have
the coordinates $(q_e,0)$, where $q_e$ satisfies the system
(\ref{sys2}). The Hessian matrix of the function $\textrm{\Large
V}$ in (\ref{V''''}) has  at the equilibrium point the same
expression (\ref{hessiana'}). By the Asymptotic Stability Theorem,
the equilibrium points are locally asymptotically stable.$\square$

\section{Example}

We consider an electrical circuit with the  associated oriented
connected graph as in Figure 1.

\begin{center}
 \scalebox{0.40}{\includegraphics*{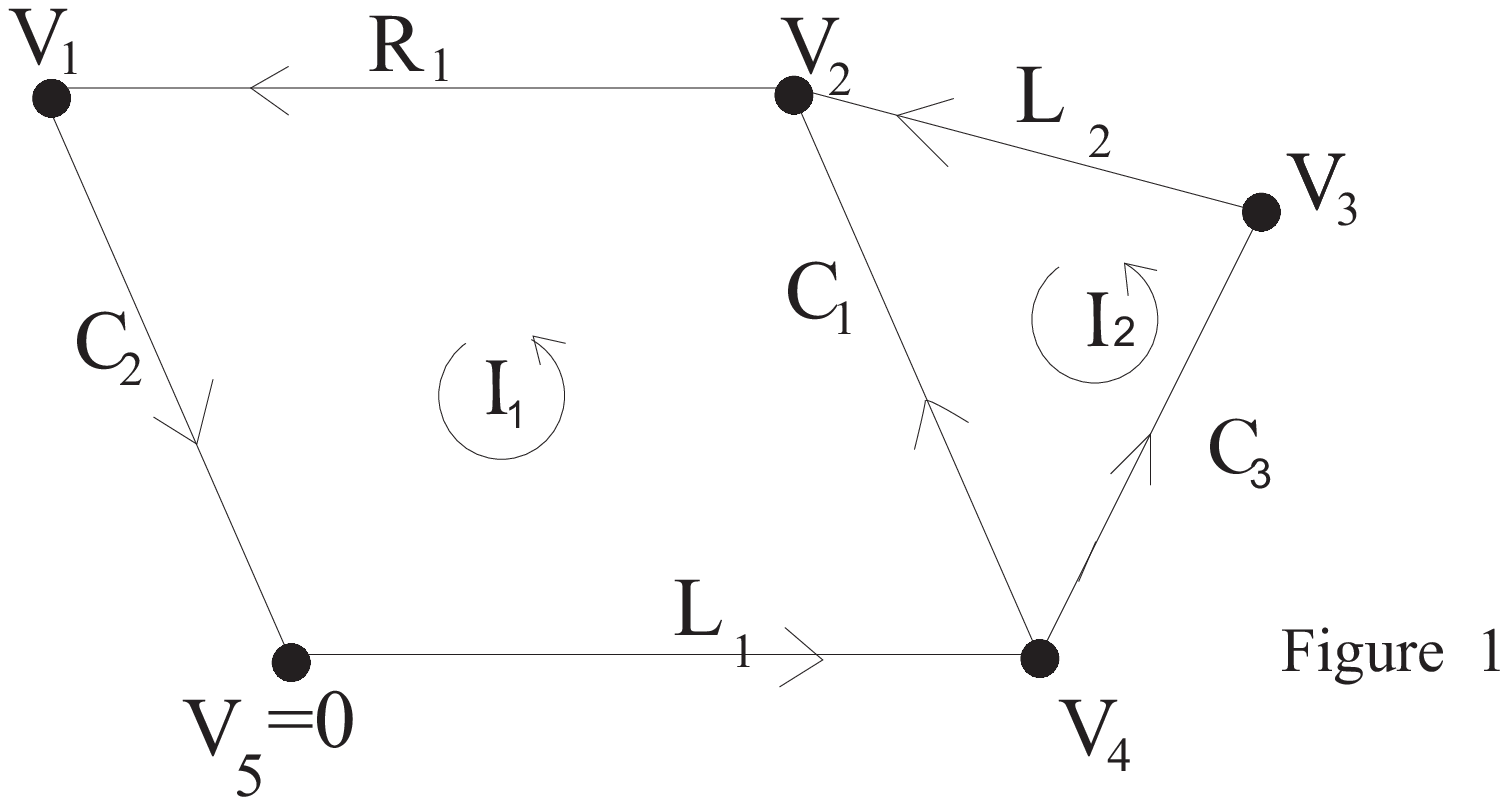}}
\end{center}

\noindent We have $r=1$, $k=2$, $p=3$, $n=4$, $m=2$ , $b=6$. We
choose the reference node to be $V_5$ and the current directions
as indicated in Figure 1. We cover the associated graph with the
loops $I_1, \, I_2$. The branches in Figure 1 are  labelled as
follows: the first branch is the resistive branch $\textsc{r}_1$,
the second and the third branch are the inductive branches
$\textsc{l}_1$, $\textsc{l}_2$ and the last three branches are the
capacitor branches $\textsc{c}_1$, $\textsc{c}_2$, $\textsc{c}_3$.
The incidence and loop matrices, $B\in
\mathfrak{M}_{64}(\mathbf{R})$ and
 $A\in \mathfrak{M}_{62}(\mathbf{R})$, write as \be
B=\left(
\begin{array}{cccc}
1&-1& 0&0\\
0& 0& 0&1\\
0& 1&-1&0 \\
0& 1& 0&-1\\
-1&0& 0&0\\
 0&0& 1& -1
\end{array}
\right),\quad \quad
 A=\left(
\begin{array}{cccc}
 1& 0\\
 1 &0 \\
0& 1  \\
 1 &-1\\
 1 &0 \\
 0 & 1
\end{array}
\right) \label{exrlc31} \ee One has  rank$(B)=4$,  rank$(A)=2$.
Kirchhoff's current law and Kirchhoff's voltage law can be
expressed by the equations \be B^T \textsc{i} =0 \quad (KCL),
\quad A^T v =0\quad (KVL)\label{kclexpl}\ee where
$\textsc{i}=(\textsc{i}_{[\Gamma]},
\textsc{i}_{(a)},\textsc{i}_{\alpha})\in \mathbf{R}\times
\mathbf{R}^2\times \mathbf{R}^3$
and $v=(v_{[\Gamma]},v_{(a)},v_{\alpha})\in \mathbf{R}\times
\mathbf{R}^2\times\mathbf{R}^3$
is the voltage drop vector.

  We
define the configuration space by
 \be M_{\mathfrak{c}}:=\{x\in \mathbf{R}^6|B^Tx=\mathfrak{c
 }\}\label{20}
\ee
 with $\mathfrak{c}$ a constant vector in $\mathbf{R}^4$.
$M_{\mathfrak{c}}$ is an affine-linear subspace in $\mathbf{R}^6$,
its dimension is $2$.
 The system in (\ref{20}) writes as
\ba
x^1-x^5=\mathfrak{c}_1\nonumber\\
-x^1+x^3+x^4=\mathfrak{c}_2\nonumber\\
-x^3+x^6=\mathfrak{c}_3\\
x^2-x^4-x^6=\mathfrak{c}_4\nonumber\label{constante} \ea We denote
local coordinates on $M_{\mathfrak{c}}$ by $q=(q^1,q^2)$. If we
take, for example,
 \be q^1:=x^5, q^2:=x^6\label{qcoord2} \ee we get
\ba
x^1&=&q^1+\mathfrak{c}_1\nonumber\\
x^2&=&q^1+\mathfrak{c}_1+\mathfrak{c}_2+\mathfrak{c}_3+\mathfrak{c}_4\nonumber\\
x^3&=&q^2-\mathfrak{c}_3\\
x^4&=&q^1-q^2+\mathfrak{c}_1+\mathfrak{c}_2+\mathfrak{c}_3\nonumber\label{exd4}
\ea Thus, the matrix of  constants $\mathcal{N}$
 in (\ref{n}) is exactly the matrix
$A$ and the constant\\
 $\mathcal{K}=\left(\begin{array}{c}
\mathfrak{c}_1\\
\mathfrak{c}_1+\mathfrak{c}_2+\mathfrak{c}_3+\mathfrak{c}_4\\
-\mathfrak{c}_3\\
\mathfrak{c}_1+\mathfrak{c}_2+\mathfrak{c}_3\\
0\\
0
\end{array}\right)$.

First we consider the case that all the electrical devices in the
circuit are linear, they are described by the relations
 (\ref{15}).
 In this case, in terms of the $q$-coordinates (\ref{qcoord2}),
the Birkhoffian $\omega_\mathfrak{c}$
on $M_\mathfrak{c}$ writes as
\ba
&&\mathfrak{Q}_1(q,\dot{q},\ddot{q})= \textrm{\scriptsize
L}_1\ddot{q}^1+\textrm{\scriptsize R}_1\dot{q}^1+
\left(\f{1}{\textrm{\scriptsize C}_1}+\f{1}{\textrm{\scriptsize
C}_2}\right)q^1-\f{1}{\textrm{\scriptsize C}_1}q^2+
\f{\mathfrak{c}_1+\mathfrak{c}_2+\mathfrak{c}_3}{\textrm{\scriptsize C}_1}\nonumber\\
&&\mathfrak{Q}_2(q,\dot{q},\ddot{q})=\textrm{\scriptsize
L}_2\ddot{q}^2 -\f{1}{\textrm{\scriptsize
C}_1}q^1+\left(\f{1}{\textrm{\scriptsize
C}_1}+\f{1}{\textrm{\scriptsize
C}_3}\right)q^2-\f{\mathfrak{c}_1+\mathfrak{c}_2+\mathfrak{c}_3}{\textrm{\scriptsize
C}_1} \label{exd6'} \ea

Let us see now how  the constants are related to the initial
conditions that may be specified for the considered
network.\\
 The differential system associated to the Birkhoffian
(\ref{exd6'}) is written \be \mathfrak{Q}_1(q,\dot{q},\ddot{q})=0,
\quad \mathfrak{Q}_2(q,\dot{q},\ddot{q})=0\label{difsist} \ee
For each capacitor we are able to specify the initial charge, that
is, $\textsc{q}_1(0),$ $\textsc{q}_2(0),$ $\textsc{q}_3(0),$ and
for each inductor the initial current, that is,
$\textsc{i}_{(1)}(0)$, $\textsc{i}_{(2)}(0)$. Taking into account
(\ref{qcoord2}), the relation $\textsc{i}=\dot{x}$ and the
relations two and three in (\ref{exd4}), we have the following
initial conditions for the differential system  (\ref{difsist})
\ba
q^1(0)&=&\textsc{q}_2(0)\nonumber\\
q^2(0)&=&\textsc{q}_3(0)\nonumber\\
\dot{q}^1(0)&=&\textsc{i}_{(1)}(0)\nonumber\\
\dot{q}^2(0)&=&\textsc{i}_{(2)}(0) \label{initial} \ea

\noindent Besides, taking into account the notations
(\ref{qcoord2}) and the last relation in (\ref{exd4}), we find \be
\mathfrak{c}_1+\mathfrak{c}_2+\mathfrak{c}_3=\textsc{q}_1(0)-\textsc{q}_2(0)+\textsc{q}_3(0)
\ee Thus, the Birkhoffian (\ref{exd6'}) becomes \ba
&&\mathfrak{Q}_1(q,\dot{q},\ddot{q})= \textrm{\scriptsize
L}_1\ddot{q}^1+\textrm{\scriptsize R}_1\dot{q}^1+
\left(\f{1}{\textrm{\scriptsize C}_1}+\f{1}{\textrm{\scriptsize C}_2}\right)q^1-\f{1}{\textrm{\scriptsize C}_1}q^2+\f{\textsc{q}_1(0)-\textsc{q}_2(0)+\textsc{q}_3(0)}{\textrm{\scriptsize C}_1}\nonumber\\
&&\mathfrak{Q}_2(q,\dot{q},\ddot{q})= \textrm{\scriptsize
L}_2\ddot{q}^2 -\f{1}{\textrm{\scriptsize
C}_1}q^1+\left(\f{1}{\textrm{\scriptsize
C}_1}+\f{1}{\textrm{\scriptsize
C}_3}\right)q^2-\f{\textsc{q}_1(0)-\textsc{q}_2(0)+
\textsc{q}_3(0)}{\textrm{\scriptsize C}_1} \label{exd6} \ea

If the constant $\textrm{\scriptsize R}_1>0$, the Birkhoffian
(\ref{exd6}) is  \textbf{dissipative}. The function
$E_{0_\omega}:TM_\mathfrak{c}\rightarrow \mathbf{R}$ and the
dissipative 1-form satisfying (\ref{dissip})  have the expressions
\ba \hspace{-0.5cm}E_{0_\omega}(q,\dot{q})&=&
\f{1}{2}\textrm{\scriptsize L}_1(\dot{q}^1)^2+
\f{1}{2}\textrm{\scriptsize L}_2(\dot{q}^2)^2+
\f{1}{2\textrm{\scriptsize C}_1}(q^1-q^2)^2+
\f{1}{2\textrm{\scriptsize C}_2}(q^1)^2+
\f{1}{2\textrm{\scriptsize C}_3}(q^2)^2+\nonumber\\
&&\f{\textsc{q}_1(0)-\textsc{q}_2(0)+\textsc{q}_3(0)}{\textrm{\scriptsize
C}_1}q^1-\f{\textsc{q}_1(0)-
\textsc{q}_2(0)+\textsc{q}_3(0)}{\textrm{\scriptsize C}_1}q^2
\label{exd7} \ea
 \ba \hspace{0.5cm} D=\textrm{\scriptsize R}_1dq^1\label{d1} \ea
The equilibrium point of the considered linear network  is
 the solution of the system
\ba &&\left(\f{1}{\textrm{\scriptsize
C}_1}+\f{1}{\textrm{\scriptsize C}_2}\right)q^1
-\f{1}{\textrm{\scriptsize
C}_1}q^2+\f{\textsc{q}_1(0)-\textsc{q}_2(0)
+\textsc{q}_3(0)}{\textrm{\scriptsize C}_1}
=0\nonumber\\
&&-\f{1}{\textrm{\scriptsize
C}_1}q^1+\left(\f{1}{\textrm{\scriptsize C}_1}+
\f{1}{\textrm{\scriptsize
C}_3}\right)q^2-\f{\textsc{q}_1(0)-\textsc{q}_2(0)+
\textsc{q}_3(0)}{\textrm{\scriptsize C}_1} =0 \label{exd8} \ea  If
the constants
 $\textrm{\scriptsize L}_1$,  $\textrm{\scriptsize L}_2$,
  $\textrm{\scriptsize C}_1$,  $\textrm{\scriptsize C}_2$,
   $\textrm{\scriptsize C}_3$ satisfy the conditions
 (\ref{condliniar'}), this equilibrium point is
  \textbf{asymptotically stable}. We define a
Liapunov function $V$ by \ba \hspace{-0.9cm}
V(q,\dot{q})&=&E_{0_\omega}(q,\dot{q})-E_{0_\omega}(q_e,0) \, =\,
\f{1}{2}\textrm{\scriptsize L}_1(\dot{q}^1)^2+
\f{1}{2}\textrm{\scriptsize L}_2(\dot{q}^2)^2\nonumber\\
&& + \, \f{1}{2\textrm{\scriptsize
C}_1}\left[(q^1-q^2)-(q^1_e-q^2_e) \right]^2+
\f{1}{2\textrm{\scriptsize C}_2}(q^1-q^1_e)^2+
\f{1}{2\textrm{\scriptsize C}_3}(q^2-q^2_e)^2\ea where $q_e$
satisfies the system (\ref{exd8}). The level curves of the
function (\ref{exd8}) represent a set of ellipsoids surrounding
the equilibrium point. Because the Birkhoffian (\ref{exd6}) is
dissipative, it follows that  $dE_{0_\omega}<0$, and therefore
$dV<0$.

Let us consider now the case that all the devices are  nonlinear,
they are described by the relations (\ref{nonliniar}).
 For the coordinate system on $M_\mathfrak{c}$ given by
(\ref{qcoord2}), the Birkhoffian becomes \ba
&&\mathfrak{Q}_1(q,\dot{q},\ddot{q})=
 L_1(\dot{q}^1)\ddot{q}^1+R_1(\dot{q}^1)+
C_1(q^1-q^2+\mathcal{K}^3)+
C_2(q^1)\nonumber\\
&&\mathfrak{Q}_2(q,\dot{q},\ddot{q})=
 L_2(\dot{q}^2)\ddot{q}^2
- C_1(q^1-q^2+\mathcal{K}^3)+ C_3(q^2) \label{exd10} \ea with
$\mathcal{K}^3=\mathfrak{c}_1+\mathfrak{c}_2+\mathfrak{c}_3=
\textsc{q}_1(0)-\textsc{q}_2(0)+\textsc{q}_3(0)$.

If $R_1$ satisfies the condition (\ref{assumption}), the
Birkhoffian (\ref{exd10}) is \textbf{dissipative}. The function
$E_{0_\omega}:TM_\mathfrak{c}\rightarrow \mathbf{R}$ and the
dissipative 1-form satisfying (\ref{dissip}) are given by
 \ba
\hspace{-0.9cm} E_{0_\omega}(q,\dot{q})&=& \int
L_1(\dot{q}^1)\dot{q}^1d\dot{q}^1+
\int L_2(\dot{q}^2)\dot{q}^2d\dot{q}^2+\int C_1(q^1-q^2+\mathcal{K}^3)(dq^1-dq^2)+\nonumber\\
&& \int C_2(q^1)dq^1+\int C_3(q^2)dq^2- \int\int
C'_1(q^1-q^2+\mathcal{K}^3)dq^1dq^2 \label{exd11} \ea  \ba
D=R_1(\dot{q}^1)dq^1\label{d1'} \ea The equilibrium points of the
considered nonlinear network are
 the solutions of the system
\ba R_1(0)+C_1(q^1-q^2+\mathcal{K}^3)+
C_2(q^1)&=&0\nonumber\\
-C_1(q^1-q^2+\mathcal{K}^3)+ C_3(q^2)&=&0 \label{exd12} \ea

 \noindent 1) If $ R_{1}(0)=0, $ and   $L_1,\,
L_2,$  $C_1,\,C_2,\,C_3$ satisfies (\ref{cond'}), then the
equilibrium points are \textbf{locally asymptotically stable}. We
define a Liapunov function $V$ by \be
V(q,\dot{q})=E_{0_\omega}(q,\dot{q})-E_{0_\omega}(q_e,0)\label{3}\ee
with $E_{0_\omega}$  given by (\ref{exd11}) and $q_e$ satisfying
(\ref{exd12}). The  Hessian matrix of $V$ at a equilibrium point
$(q_e,0)$ has the expression \be \textbf{H}_{V}(q_e,0)=
\left(\begin{array}{cccc}
L_1(0) & 0 &0 &0 \\
0 & L_2(0)& 0 &0 \\
0&0&\tilde{C}'_1(q_e^1,q_e^2)+
C'_2(q_e^2)&-\tilde{C}'_1(q_e^1,q_e^2)\\
0&0&-\tilde{C}'_1(q_e^1,q_e^2)& \tilde{C}'_1(q_e^1,q_e^2)+
C'_3(q_e^2)
\end{array}
\right) \label{exd13} \ee Under the assumptions we made, this
matrix is positive definite. The centers of the level curves of
the function (\ref{3}) have the coordinates $(q_e,0)$, where $q_e$
satisfies the system (\ref{exd12}) with $R_1(0)=0$. Because the
Birkhoffian (\ref{exd10}) is dissipative, it follows that
$dE_{0_\omega}<0$, and therefore $dV<0$.

\noindent 2) If  $ R_{1}(0)\neq 0, $ but
 for all $x\neq 0$ \be x\left(R_{1}(x)-R_{1}(0)\right)
>0 \label{second'}\ee
and  $L_1,\, L_2,$ $C_1,\,C_2,\,C_3$ satisfies (\ref{cond'}), then
the equilibrium points are \textbf{locally asymptotically stable}.
We define now a Liapunov function  $\textrm{\Large V}$ by \be
\textrm{\Large V}(q,\dot{q})=\textrm{\Large
E}_{0_\omega}(q,\dot{q})-\textrm{\Large
E}_{0_\omega}(q_e,0)\label{energienonlinear''''}\ee with {\Large
E}$_{0_\omega}:TM_\mathfrak{c} \rightarrow \mathbf{R}$ given by
\be \textrm{\Large
E}_{0_\omega}(q,\dot{q})=E_{0_\omega}(q,\dot{q})+R_{1}(0)q^1
\label{energienonlinear'''} \ee $E_{0_\omega}(q,\dot{q})$  having
the expression (\ref{exd11}). The centers of the level curves of
the function (\ref{energienonlinear''''}) have the coordinates
$(q_e,0)$, where $q_e$ satisfies the system (\ref{exd12}). The
Hessian matrix of the function $\textrm{\Large V}$ in
(\ref{energienonlinear''''}) has  at the equilibrium point the
same expression (\ref{exd13}). It remains to prove that
$d\textrm{\Large V}<0$. This yields from the dissipativeness of
the Birkhoffian (\ref{exd10}). We consider the following
dissipative 1-form \be \textrm{\Large
D}=\left[R_1\left(\dot{q}^1\right)-R_{1}(0)\right]dq^1\label{seconddissip'}\ee
In  view of assumption (\ref{second'}), this vertical 1-form is
indeed dissipative.
 One can easily check that  for the
function $\textrm{\Large E}_{0_\omega}(q,\dot{q})$ in
(\ref{energienonlinear'''}) and the dissipative 1-form in
(\ref{seconddissip'}), the following identity is fulfilled \be
\sum^{2}_{j=1}Q_j(q\, \, \dot{q},\,
\ddot{q})\dot{q}^j=\sum^{2}_{j=1}\left[\f{\pa \textrm{\Large
E}_{0_\omega}}{\pa q^j}\dot{q}^j+\f{\pa \textrm{\Large
E}_{0_\omega}}{\pa \dot{q}^j}\ddot{q}^j+ \textrm{\Large
D}_{j}(q,\dot{q})\dot{q}^j \right] \ee
  The Birkhoffian (\ref{exd10})
being dissipative we have $d\textrm{\Large E}_{0_\omega}<0$,
therefore $d\textrm{\Large V}<0$.


\textit{\textbf {Acknowledgement.}} This work  was  supported by
IMAR through the contract of excellency CEx  06-11-12/ 25.07.06.

\end{document}